\documentclass[12pt]{article}
\usepackage{multicol,amsmath,amssymb,amstext,graphicx}
\newcommand{\no}[1]{\|{#1}\|}

\def\nn{\nonumber}

\newcommand{\h}{\hspace*{.24in}}

\def\bq{\begin{equation}}
\def\eq{\end{equation}}

\newcommand{\ep}{\epsilon}

\newcommand{\be}{\begin{eqnarray*}}
\newcommand{\en}{\end{eqnarray*}}
\newcommand{\bes}{\begin{eqnarray}}
\newcommand{\ens}{\end{eqnarray}}

\newcommand{\et}{ e^{(p-1)T \xi^2}}
\newcommand{\eg}{ e^{pT \xi^2}}
\newcommand{\ega}{ \beta e^{pT \xi^2}}
\newcommand{\ef}{ \exp\{(p-1)T \xi^2\}}

\begin{document}
\title{\bf \LARGE  A new stability results for the backward heat equation   } %
\author{  {\bf Alain Pham Ngoc Dinh}\\MAPMO, UMR CNRS 6628, Orleans University, France.
\and{\bf Dang Duc Trong}\\Department of Mathematics ,HoChiMinh City
National University,  \\227 Nguyen Van Cu, Q. 5, HoChiMinh City,
VietNam. \and {\bf Pham Hoang Quan}\\ Department of Mathematics, Sai
Gon University, Ho Chi Minh city, Viet Nam \and{\bf Nguyen Huy Tuan}
\\Department of Mathematics and Informatics,Ton Duc Thang
University, \\ 98, Ngo Tat To, Binh Thanh district, Ho Chi Minh city
, Viet Nam
 }
\maketitle%

\begin{center}{\bf Abstract} \end{center}\par
{In this paper, we regularize  the nonlinear  inverse time  heat   problem
in the unbounded region by Fourier method. Some new convergence rates are obtained.
Meanwhile, some quite sharp error estimates between the approximate
solution and exact solution are provided. Especially, the optimal convergence of
the approximate solution at t = 0 is also proved. This  work extends to many  earlier results in \cite{f2,f3, hao1,Quan,tau1, tau2,  Trong3,x1}.\\
{\bf Keywords and phrases:}  Backward  heat  problem,  Ill-posed problem, Fourier Transform, Contraction principle.

{\bf Mathematics subject Classification 2000:} 35K05, 35K99, 47J06, 47H10.\\

\section{ Introduction}  \par
Transient heat conduction phenomena are generally governed by the parabolic heat conduction
equation and if the initial temperature distribution and boundary conditions are
specified, then this, in general, leads to a well posed problem which may easily be solved
numerically by using various  methods.

However, in many practical situations when dealing with a heat conducting body it is not
always possible to specify the boundary conditions or the initial temperature. For example, in
practice, one may have to investigate the temperature distribution and the heat flux history
from the known data at a particular time. In other words, it may be possible to specify the
temperature distribution at a particular time, say $t = T > 0,$  and from this data the question arises as
to whether the temperature distribution at any earlier time $t < T$ can be retrieved. This is
usually referred to as the backward heat conduction problem (BHCP), or the final boundary
value problem.
The systematic study of the backward heat conduction problem is of rather recent origin,
although isolated considerations have been given to such problems for several hundred years.
In general no solution which satisfies the heat conduction equation, the final data and the
boundary conditions exists. Further, even if a solution did exist, it would not be continuously
dependent on the boundary and the final data, see Payne [2]. Thus the BHCP is an
example of an ill-posed problem that is impossible to solve using classical numerical methods
and requires special techniques to be employed, see Hadamard [3]. Conditions for which the
BHCP becomes well-posed have been investigated by Miranker [4] and Cannon and Douglas
[5]. These studies introduced additional hypotheses which restrict the class of functions to
which the solution must belong, and which are seldom satisfied. Therefore numerical methods
of solution appear more useful. Thus regularization techniques, see for example Cannon [6]
and Han et al. [7], have been developed for solving the BHCP. Different methods, based on
a perturbation of the original parabolic heat equation were proposed by Lattes and Lions [8]
and Lesnic et al. [9]. These methods consist of replacing the operator by a perturbed higher
order one that holds better invertibility features.

In this paper, we   consider the following  problem
\bq
\left\{ \begin{gathered}
  u_t  - u_{xx}  = 0,~~~(x,t) \in R \times (0,T), \hfill \\
  u(x,T) =\varphi(x),~~x\in R,\hfill\\
 \end{gathered}  \right.\label{1}
\eq
where  $T >0$ be a given and  $\varphi(x)$ are given.
We want to retrieve the temperature distribution $u(x, t)$ for $0 \le  t < T.$ Of course, since the data $\varphi(.)$ is based on
(physical) observations, there will be measurement errors, and we would actually have as data some function $\varphi_\ep \in L^2(R)$
, for which
$\|\varphi-\varphi_\ep\| \le \ep$, where $\|.\|$ denotes the $L^2$-norm, the constant $\ep>0$ represents a bound on the measurement error. That is to
say, practically, we need to consider the following problem,
\bq
\left\{ \begin{gathered}
  u_t  - u_{xx}  = 0,~~~(x,t) \in R \times (0,T), \hfill \\
  u(x,T) =\varphi_\ep(x),~~x\in R,\hfill\\
 \end{gathered}  \right.\label{2}
\eq
 Notice the reader that the  problem (1) is investigated  in some recent papers of ChuLiFu \cite{f2,f3,x1} and of other authors such as Lien\cite{Trong3},Murniz \cite{m3}, et al . To the authors' knowledge, so far there are many papers on the backward heat equation, but
theoretically the error estimates of most regularization methods in the literature are Holder type,
i.e.,
\bes
\|u(.,t)-v^\ep(.,t)\| \le C \ep^{k},~~~k>0 \label{ht18}.
\ens
where $C$ is the constant depend on  $u$,  $k$ is a constant is not depend on $t, u$ . As we know, $\ep^{k}$ converges to zero more quickly than the logarithmic term. So,
the major object of this paper is to provide  new
regularization method to established the Holder estimates such as \eqref{ht18}. We give   a new approximation  problem  and investigate the
error estimate between the regularization solution and the exact one. \\
The remainder of the paper is divided into three sections. In Section 2, we establish
the approximated problem and show that it is well posed. Then, we also estimate the
error between an exact solution u of Problem (1) and the approximation solution $ v^\ep(.,t)$ with the Holder type. Finally, a numerical experiment will be given in Section 3.
\section{Regularization and error estimates.}
Let
\be
\hat \varphi(\xi)=\frac{1}{\sqrt{2\pi}} \int_{-\infty}^{+\infty}\varphi(x)e^{-i\xi x}dx
\en
be the Fourier transform of the  function $\varphi \in L^2(R)$.
We assume $u(x, t)$ is the unique
solution of \eqref{1}. Using the Fourier transform technique to problem \eqref{1} with respect to the
variable $x$, we can get the Fourier transform  $\hat u(\xi,t) $ of the exact solution $u(x, t) $ of problem \eqref{1}
\bq
\left\{ \begin{gathered}
  \hat u_t(\xi,t)= (i\xi)^2u(\xi,t), \hfill \\
\hat  u(\xi,T) =\hat \varphi(\xi),~~\xi \in R,\hfill\\
 \end{gathered}  \right.\label{3}
\eq
The solution to equation \eqref{3} is given by
\bes
\h\h\h\h\h\hat u(\xi,t) = e^{(T-t)\xi^2   } \hat \varphi (\xi) . \label{q33}
\ens
or equivalently,
\bes
u(x,t)& =& \frac {1}{\sqrt{2 \pi}} \int \limits_{-\infty}^{+\infty} {  e^{(T-t)\xi^2   }   \hat \varphi(\xi) e^{i\xi x}d\xi}.\label{qqqq4}
\ens
Since $t < T,$ we know from \eqref{qqqq4} that, when $\xi$ becomes large, $\exp{(T -t)\xi^2} $      increases rather quickly. Thus
for $   \hat u(\xi,t) \in L^2(R) $    with respect to $\xi$, the exact data function $  \hat \varphi(\xi) $ must decay rapidly as $ |\xi| \to \infty$. Small errors in high-frequency components can blow up and completely destroy the solution for $0 \le  t < T$.As for the
measured data $\varphi_\ep(x)$, its Fourier transform $\hat \varphi_\ep(\xi ) $ is merely in $ L^2(R)$.

Let $p>1$ be a constant number.
We approximated problem (3) by pertubing the Fourier transform of  final value $\varphi$ as follows
 \bq
\left\{ \begin{gathered}
  \frac{\partial  \hat w_\ep(\xi,t)}{\partial t}= (i\xi)^2\hat w_\ep(\xi,t), \hfill \\
\hat  w_\ep(\xi,T) = \frac{e^{-T\xi^2}}{\beta \et+e^{-T\xi^2}}\hat \varphi (\xi) ,~~\xi \in R,\hfill\\
 \end{gathered}  \right.\label{4}
\eq
The formal solution of \eqref{4} is also easily seen to be
\bes
\h\h\h\h\hat w_\ep(\xi,t) = \frac{e^{-t\xi^2}}{\beta \et+e^{-T\xi^2}}\hat \varphi (\xi),\label{qqq3}
\ens
or
\bes
w_\ep(x,t)& =& \frac {1}{\sqrt{2 \pi}} \int \limits_{-\infty}^{+\infty} {\frac{e^{-t\xi^2}}{\beta \et + e^{-T\xi^2}}\hat \varphi(\xi) e^{i\xi x}d\xi} ,\label{q44}
\ens
where $\beta$ is a positive number such that $ \beta>0 $.\\
Let $v_\ep(.,t)$ be the approximated solution given by
\bes
v_\ep(x,t)& =& \frac {1}{\sqrt{2 \pi}} \int \limits_{-\infty}^{+\infty} {\frac{e^{-t\xi^2}}{\beta \et + e^{-T\xi^2}}\hat \varphi_\ep(\xi) e^{i\xi x}d\xi} ,\label{q45}
\ens
Note that if $\beta$ is chosen small, then for small $|\xi|$, $ \frac{e^{-t\xi^2}}{\beta \et+e^{-T\xi^2}}  $
in \eqref{qqq3} is close to $  e^{(T-t)\xi^2   }    $ in \eqref{q33}.

For $\xi,x,\beta >0, 0\le a\le b$, we  prove the following inequality
\bes
&&i) .  \frac{e^{a\xi^2}}{1+\beta  e^{b\xi^2}} \le \beta^{-\frac{a}{b}}.
\label{ht4}\\
&&ii).  \frac{e^{a\xi^2}}{\xi^2(1+\beta  e^{b\xi^2})} \le \frac{b}{\ln(\frac{1}{\beta})}\beta^{-\frac{a}{b}} \label{ht5}
\ens
Thus, we have
\be
\frac{e^{a\xi^2}}{1+\beta  e^{b\xi^2}} &=&\frac{e^{a\xi^2}}{(1+\beta  e^{b\xi^2})^{\frac{a}{b}}  (1+\beta  e^{b\xi^2})^{1-\frac{a}{b}}         } \\
&\le& \frac{e^{a\xi^2}}{(1+\beta  e^{b\xi^2})^{\frac{a}{b}} }\\
& \le& \beta ^{-\frac{a}{b}}.
\en
For $0 \le t \le s \le T$, denote
\be
A(\xi,t)&=&  \frac{\exp\{-t\xi^2\}}{\beta \ef+ \exp\{-T\xi^2\}},\\
B(\xi,s,t)&=&\frac{\exp\{(s-t-T)\xi^2\}}{\beta \ef+\exp\{-T\xi^2\}}.\\
\en
Using the inequality \eqref{ht4}, we obtain
\bes
A(\xi, t) &=& \frac {{e^{(T-t) \xi^2}}}{1+\beta e^{p T\xi^2}}\le \beta^{\frac   {t-T}{pT}     }
\label {m1} \\
B(\xi, s,t) &=&  \frac {{e^{(s-t) \xi^2}}}{1+\beta e^{p T\xi^2}}\le \beta^{\frac   {t-s}{pT}     } \label{m2}.
\ens

{\bf Theorem 1.}
\emph{ Let $\varphi\in L^2(R)$. Then unique solution of the problem \eqref{qqqq4}   depends on the final value $\varphi$,i.e, if $w,v$ are the solution of  the  problem \eqref{qqqq4}corresponding to the final values $\varphi $ and $\omega$, then
\begin{eqnarray*}
\no{w(.,t)-v(.,t)}\le  \beta^{\frac   {t-T}{pT}     }\no{\varphi-  \phi}.
\end{eqnarray*}
}
{\bf Proof of Theorem  1.}\\
Let $w$ and $v$ be two solution of the problem \eqref{qqqq4} corresponding to the final values $\varphi$ and $ \phi$. Using Parseval inequality, we have
\be
\|w(.,t)-v(.,t)\|^2 &=&\|\hat w(.,t)-\hat v(.,t)\|^2\\
&\le&  \int \limits_{- \infty}^{+\infty} \left|A(\xi,t)\left( \hat \varphi(\xi)-  \hat \phi(\xi)            \right)\right|^2d \xi\\
 &\le &\beta^{\frac   {2t-2T}{pT}     } \int \limits_{- \infty}^{+\infty} \left|\left( \hat \varphi(\xi)-  \hat \phi(\xi)            \right)\right|^2d \xi\\
 &\le& \beta^{\frac   {2t-2T}{pT}     }   \|\hat \varphi-  \hat \phi  \|^2\\
&\le&\beta^{\frac   {2t-2T}{pT}     } \| \varphi-   \phi  \|^2.
\en
Hence
\begin{eqnarray*}
\no{w(.,t)-v(.,t)}\le  \beta^{\frac   {t-T}{pT}     }\no{\varphi-  \phi}.
\end{eqnarray*}
{\bf Remark 1.}\\
In \cite{clark, Trong1}, the stability of magnitude is $e^{\frac{T}{\epsilon}}$.
In \cite {f2,Quan, Trong2}, the stability estimate is of order $ \epsilon^{\frac{t}{T}-1}$.\\
In our paper, we give a better estimation of the stability order, which is $C\beta^{\frac{t}{T}-1}\left(\frac{T}{ 1+\ln(\frac{T}{\beta})}\right)^{1-\frac{t}{T}}$. It is easy to see that the order of the error, introduced by small changes in the final value $g$, is less than the order given in \cite{Quan, Trong2}. This is among of  the advantages of our method.\\
{\bf Theorem 2.}
{\it  Let $w_\ep$ and $v_\ep$ defined by \eqref{q44} and \eqref{q45}. Then one has
\bes
\|w_\ep(.,t)-v_\ep(.,t)\| \le \beta^{\frac   {t-T}{pT}     } \ep.
\ens
}
{\bf Proof of Theorem 2.}
Using Theorem 1, we get
\bes
\|w_\ep(.,t)-v_\ep(.,t)\| \le \beta^{\frac   {t-T}{pT}     } \no{\varphi_\ep-  \varphi}\le \beta^{\frac   {t-T}{pT}     } \ep.
\ens
{\bf Theorem 3.}
{\it Let $u$ be the exact solution of problem such that $\|u(.,0)\| \le E_1 $. Then one has
\be
\|u(.,t)-v_{\ep } (.,t)\| \le  \ep^{\frac{t}{T}} ( E_1+1).
\en
}
{\bf Proof of Theorem 2.}
Using \eqref{ht4}, we have
\begin{eqnarray*}
\|u(.,t)-w_{\ep } (.,t)\|^2
&=& \int \limits_{-\infty}^{+\infty}{|\hat u(\xi,t)-\hat w_\beta (\xi,t)|^2d\xi}\\
&=&\int \limits_{-\infty}^{+\infty}  {\left| {(e^{(T-t)\xi^2}-A(\xi,t))\hat \varphi (\xi)
 }     \right|^2 }d\xi  \\
&=&\int \limits_{-\infty}^{+\infty}{\left| {\frac{\beta \eg e^{(T-t)\xi^2}}{(1 + \ega)}\hat \varphi (\xi)
} \right.}  \\
&=&\int \limits_{-\infty}^{+\infty} \left| {\frac{\beta e^{(pT-t)\xi^2} }{(1+\ega)}\hat u (\xi,0)
}\right|^2 d\xi    \\
&\le & \beta^2 \beta^{\frac   {2t-2T}{pT}     }  \int \limits_{-\infty}^{+\infty} \left| {\hat u (\xi,0)
}\right|^2 d\xi \\
&\le& \beta^{\frac   {2t}{pT}     }E_1^2.
\en
Therefore
\be
\|u(.,t)-w_{\ep} (.,t)\| \le  \beta^{\frac   {t}{pT}     }E_1.
\en
From $\beta=\ep^{p}$ and using the inequality, we obtain
\be
\|u(.,t)-v_{\ep} (.,t)\| &\le& \|u(.,t)-w_{\ep} (.,t)\| +\|w_\ep(.,t)-v_{\ep} (.,t)\| \\
&\le& \beta^{\frac   {t}{pT}     }E_1+\beta^{\frac   {t-T}{pT}     } \ep\\
&\le & \ep^{\frac{t}{T}} ( E_1+1).
\en
{\bf Remark 2. }
Notice that the convergence estimate in Theorem 1 does not give any useful information on the continuous
dependence of the solution at $t = 0.$ This is common in the theory of ill-posed problems, if we do not
have additional conditions on the smoothness of the solution. To retain the continuous dependence of the
solution at $ t = 0,$ instead of (2.5), one has to introduce a stronger a priori assumption.

We denote  $\|.\|_k$ be  the norm in Sobolev space  $ H^k(R),~k>0$ defined by
\be
\|u(.,0)\|_{k}:=\left( \int_{-\infty}^{+\infty} (1+\xi^2)^k   |\hat u (\xi,0)|^2 d\xi \right)^{\frac{1}{2}}.
\en
{\bf Theorem 3.}
{\it Let $u$ be the exact solution of problem such that $ \|u(.,0)\|_{2} \le E_2   $. Let $\beta=\ep$ , then one has
\bes
\|u(.,t)-w_{\ep } (.,t)\| \le  \frac{pT}{\ln(\frac{1}{\ep})}\ep^{\frac{t}{pT}} E_2+\ep^{\frac   {t-T+pT}{pT}     }. \label{qtt1}
\ens
}
{\bf Proof of Theorem 3.}\\
\be
\|u(.,t)-w_{\ep } (.,t)\|^2 &=&\int \limits_{-\infty}^{+\infty} \left| {\frac{\beta e^{(pT-t)\xi^2} }{\xi^2(1+\ega)}\xi^2\hat u (\xi,0)
}\right|^2 d\xi   \\
&\le&    \int \limits_{-\infty}^{+\infty} \xi^4\hat u ^2(\xi,0) d\xi\\
&\le &\left(\frac{pT\beta}{\ln(\frac{1}{\beta})}\beta^{\frac{t-pT}{pT}}\right)^2 \int \limits_{-\infty}^{+\infty} \xi^4\hat u ^2(\xi,0) d\xi\\
&\le& \left(\frac{pT\beta}{\ln(\frac{1}{\beta})}\beta^{\frac{t-pT}{pT}}\right)^2 \|u(.,0)\|_{2}^2\\
&\le &\frac{pT}{\ln(\frac{1}{\beta})}\beta^{\frac{t}{pT}} E_2.
\en
From $\beta=\ep$ and using the inequality, we obtain
\be
\|u(.,t)-v_{\ep} (.,t)\|  &\le& \|u(.,t)-w_{\ep} (.,t)\| +\|w_\ep(.,t)-v_{\ep} (.,t)\| \\
&\le& \frac{pT}{\ln(\frac{1}{\ep})}\ep^{\frac{t}{pT}} E_2+\ep^{\frac   {t-T+pT}{pT}     } .
\en
{\bf Remark 3. }\\
1. In $t=0$, the error \eqref{qtt1} becomes
\bes
\|u(.,0)-v_{\ep} (.,0)\|
&\le& \frac{pT}{\ln(\frac{1}{\ep})} E_2+\ep^{\frac   {-T+pT}{pT}     } \label{qtt4} .
\ens
It follows from  $p>1$ that the right hand side of \eqref{qtt4} converges to zero when $\ep \to 0$. This error is the same order as Theorem 3.1 in paper \cite{f3} (see page 566) and Theorem 2.1 in paper \cite{f4}.

2.Since \eqref{qtt1}, the first term of the right hand side of \eqref{qtt1} is the logarithmic form, and the second term is a power, so the order of \eqref{qtt1} is also logarithmic order. It is the same order as some results which is mentioned in Remark 2. This often occurs in the boundary error estimate
for ill-posed problems. To retain the Holder order in $[0,T]$,  we
introduce a different  priori assumption.\\
{\bf Theorem 4.}\\
{\it Assume that there exist a positive number $\gamma \in (0,pT)$ such that
\be
\int \limits_{-\infty}^{+\infty} e^{2\gamma \xi^2} \hat u (\xi,0)d\xi <E_3^2.
\en
Let $\beta=\ep$ and $h=\min\{ \gamma, (p-1)T\}$. Then, one has
\bes
\|u(.,t)-w_{\ep } (.,t)\| \le  \ep^{\frac{t+h}{pT}} (E_3+1).  \label{t22}
\ens
}
{\bf Proof of Theorem 4.}\\
From , we get
\bes
\|u(.,t)-w_{\ep } (.,t)\|^2
&=& \int \limits_{-\infty}^{+\infty}{|\hat u(\xi,t)-\hat w_\beta (\xi,t)|^2d\xi}\nn\\
&=&\int \limits_{-\infty}^{+\infty}{\left| {\frac{\beta \eg e^{(T-t)\xi^2}}{(1 + \ega)}\hat \varphi (\xi)
} \right.}  \nn\\
&=&\int \limits_{-\infty}^{+\infty} \left| {\frac{\beta \eg }{1+\ega}\hat u (\xi,t)
}\right|^2 d\xi    \nn\\
&=&\int \limits_{-\infty}^{+\infty} \left| {\frac{\beta e^{(pT-t-\gamma)\xi^2} }{1+\ega}e^{\gamma \xi^2}\hat u (\xi,0)
}\right|^2 d\xi. \label{ht7}
\ens
Using the inequality $\eqref{ht4}$, we have
\bes
\frac{\beta e^{(pT-t-\gamma)\xi^2} }{1+\ega} \le \beta \beta^{\frac{t+\gamma-pT}{pT}} .\label{ht8}
\ens
Combining  \eqref{ht7} and \eqref{ht8}, we get
\be
\|u(.,t)-w_{\ep } (.,t)\|^2  \le \beta^{\frac{2t+2\gamma}{pT}}  \int \limits_{-\infty}^{+\infty} e^{2\gamma \xi^2} \hat u (\xi,0)d\xi \le \beta^{\frac{2t+2\gamma}{pT}} E_3^2.
\en
From $\beta=\ep$ and using the inequality, we obtain
\be
\|u(.,t)-v_{\ep} (.,t)\|  &\le& \|u(.,t)-w_{\ep} (.,t)\| +\|w_\ep(.,t)-v_{\ep} (.,t)\| \\
&\le& \beta^{\frac{t+\gamma}{pT}} E_3+\beta^{\frac   {t-T}{pT}     } \ep \\
&\le &   \ep^{\frac{t+h}{pT}} (E_3+1)        .
\en
{\bf Remark 4.}\\
1. In $t=0$, the error \eqref{t22} becomes
\bes
\|u(.,t)-w_{\ep } (.,t)\| \le  \ep^{\frac{h}{pT}} (E_3+1). \label{ht11}
\ens
2. Suppose that $E_\ep=\|v^\ep-u\|$ be the error of the exact solution and the approximate solution.
In most of results concerning the backward heat, then  optimal error  between  is of the logarithmic form. It means that
\be
E_\ep \le C  \left(\ln{\frac{T}{\ep}}\right)^{-m}
\en
where $m>0$.
The error order of logarithmic form is investigated in many recent papers, such as \cite{clark, f3, Denche,  Trong1, Trong2, Trong4}.

To illustrate this, we can enumerate some more recent papers considering the errors of logarithmic order.\\
 Let $u$ and $v^\ep$ be  exact solution and the approximated solution respectively.\\
In \cite{  Trong2, Trong4} , the error is given the form
\be
\|v^\ep (.,t)-u(.,t)\| \le \frac{C_1}{1+\ln{\frac{T}{\ep}}}.
\en
In recently, Feng Xiao-Li and coauthors \cite{feng} gave the error estimates as follows
\be
\no{v^{\alpha,\delta}(.)-u(.)  } \le \frac{\delta}{2\sqrt{\alpha}}+\max\{ \left(\frac{4T}{\ln{\frac{1}{\alpha}}}\right)^{\frac{p}{2}} ,\alpha^{\frac{1}{2}} \}.
\en
 From the discussed error, we can see that most of recent regularization methods established the logarithimic stability. The convergence rates in here is very slowly.\\
In Tautenhahn and  Schr¨oter  \cite{tau2}, the authors proved that the best possible worst
case error for identifying u(t)  is given by
$
E_\ep=D \ep^{\frac{t}{T}}.
$\\
In 2007,  Zhi Qian et al. \cite{f3} (See Remark 3.6, p.570)gave the error estimation in $t=0$ for (2) as follows
\be
\|u(.,0-v(.,0\| \le \frac{1}{(\ln(1/\epsilon))^{\frac{8}{5}}}+ \max \{1,T\} \beta^{\frac{1}{5}}E. \label{hay1}
\en
where $ \beta=\frac{T}{\ln \Big(    \frac{1}{\epsilon}(\ln\frac{1}{\epsilon})^{\frac{-8}{5}}      \Big)}$.\\
In \cite{f4}, ChuLiFu and  his coauthors esrablished the logarithimic order of the form
\be
&&\|u(.,t)-u_{\delta,\xi_{\max}}\|\\
&&~~~~~~~~~~~~~\le E^{1-\frac{t}{T}} \left (\ln \frac{E}{\delta} \right)^{\frac{(t-T)s}{2T}} \left(            1+  \left ( \frac{\ln \frac{E}{\delta}  }          {  \frac{1}{T} \ln \frac{E}{\delta}+\ln(\ln \frac{E}{\delta})^{\frac{-s}{2T}}       } \right)^{\frac{s}{2}}\right).
\en
Comparing \eqref{ht11}  with the discussed error, we can see  \eqref{ht11} is the optimal error.

\section{Numerical results.}


\begin{thebibliography}{99}

\bibitem{a4} K.A. Ames, L.E. Payne. (1998),
{ Asymptovic for two regularizations of the Cauchy problem for the
backward heat equation},{\it  Math.Models Methods Appl. Sci. }
187-202.


\bibitem {h3} B. M. Campbell H  and R.J. Hughes  (2007),
{ Continuous dependence results for inhomogeneous ill-posed problems in Banach space
             }, {\it  J. Math. Anal. Appl.,} Volume 331, Issue 1,  Pages 342-357


\bibitem {h4}B. M. Campbell H, R. Hughes, E. McNabb.(2008), { Regularization of the backward heat equation via heatlets}, {\it Electron. J. Diff. Eqns.,} Vol. 2008 , No.130, 2008,
pp. 1-8.


\bibitem{clark} G. W.  Clark, S. F. Oppenheimer;
\emph{Quasireversibility methods for non-well posed problems},
Elect. J. Diff. Eqns., \textbf{1994} (1994) no. 8, 1-9.



\bibitem{Denche}  Denche,M. and Bessila,K.,   \emph {  A modified quasi-boundary value method for ill-posed problems}, J.Math.Anal.Appl,  Vol.{\bf 301}, 2005,  pp.419-426.

\bibitem{feng} Feng, Xiao-Li; Qian, Zhi; Fu, Chu-Li {\it Numerical approximation of solution of nonhomogeneous backward heat conduction problem in bounded region}.  Math. Comput. Simulation  79  (2008),  no. 2, 177--188


\bibitem{f2} Chu-Li Fu, Xiang-Tuan Xiong,and Zhi Qian. (2007),
{On  three spectral regularization method for
a backward heat conduction problem.}{\it  J. Korean Math. Soc. }{\bf 44} , No. 6, pp. 1281-1290.

\bibitem{f3} Chu-Li Fu , Zhi Qian ,  Rui Shi. (2007),{A modified method for a backward heat conduction problem}, {\it  Applied Mathematics and Computation}, {\bf 185}  564-573.

\bibitem {f4}Fu, Chu-Li; Xiong, Xiang-Tuan; Qian, Zhi {\it  Fourier regularization for a backward heat equation.}  J. Math. Anal. Appl.  331  (2007),  no. 1, 472--480


\bibitem{hao1} Dinh Nho Hao and  Nguyen Van Duc.(2009), \emph{Stability results for the heat equation backward in time},  to appear in {\it   J. Math. Anal. Appl. }



\bibitem{lattes} R. Latt\`{e}s, J.-L. Lions; \emph{M\'ethode de
Quasi-r\'eversibilit\'e et Applications}, Dunod, Paris, 1967.


\bibitem{showalter1} R.E. Showalter;
\emph{The final value problem for evolution equations}, J. Math.
Anal. Appl, \textbf{47} (1974), 563-572.










\bibitem{m2} I. V. Mel'nikova, Q. Zheng and J. Zheng;
\emph{Regularization of weakly ill-posed Cauchy problem}, J.
Inv. Ill-posed Problems, Vol. \textbf{10} (2002), No. 5, 385-393.

\bibitem {m3} Muniz W B.\emph{ A comparison of some inverse methods for estimating the initial condition of the heat
equation.} J. Comput. Appl. Math.,1999,103: 145-163


\bibitem{Quan} Quan, P. H. and Trong, D. D., \emph{ A nonlinearly backward
heat problem: uniqueness, regularization and error estimate},
Applicable Analysis, Vol. 85, Nos. 6-7, June-July 2006, pp. 641-657.


\bibitem {seidman} T.I. Seidman, \emph{Optimal filtering for the backward heat equation, }SIAM J. Numer. Anal. 33 (1996) 162–170.


\bibitem{tau1}Tautenhahn U 1998\emph{ Optimality for ill-posed problems under general source conditions} Numer. Funct. Anal.
Optim. 19  377-398.

\bibitem {tau2} U. Tautenhahn and T. Schr¨oter. (1996), On optimal regularization methods for the backward heat equation,{\it
Zeitschrift f¨ur Analysis und ihre Anwendungen.} 15 , no. 2, 475-493.


\bibitem {Trong1} D.D.Trong and N.H.Tuan \emph{ Regularization and
error estimates for nonhomogeneous backward heat problems},
Electron. J. Diff. Eqns., Vol. 2006 , No. 04, 2006,
pp. 1-10.

\bibitem {Trong2} Dang Duc Trong, Pham Hoang Quan, Tran Vu Khanh and Nguyen Huy Tuan{\it  A nonlinear case of the 1-D backward heat problem:
Regularization and error estimate}, Zeitschrift  Analysis und ihre Anwendungen, Volume 26, Issue 2, 2007, pp. 231-245.

\bibitem {Trong3}Dang Duc Trong and Tran Ngoc Lien.(2007), {Regularization of a discrete backward problem
 using coefficients of truncated Lagrange polynomials},{\it Electron. J. Diff. Eqns.,}Vol. 2007(2007), No. 51, pp. 1--14.\


\bibitem {Trong4}  Dang Duc Trong and Nguyen Huy Tuan  \emph{A nonhomogeneous backward heat problem: Regularization and error estimates }, Electron. J. Diff. Eqns., Vol. 2008 , No. 33,
pp. 1-14.

\bibitem {t5}D.D. Trong and N.H. Tuan  (2009),  { Regularization of the nonlinear backward heat problem using a method of integral equation }, {\it accepted for
publication in Nonlinear Analysis, Theory, Methods and Applications, Series A: Theory
and Methods.}\


\bibitem {x1}Xiang-Tuan Xiong, Chu-Li-Fu, Zhi Qian, and  Xiang Gao (2006),{ Error  estimates of a difference
approximation method for a backward heat conduction problem},{\it International Journal of Mathematics and Mathematical Sciences}.
Volume 2006, Article ID 45489, Pages 1-9.



\end{thebibliography}
\end{document}